\documentclass[12pt,twoside]{article}
\usepackage{amsmath}
\usepackage{amssymb}
\usepackage{enumerate,epsfig}
\usepackage[english,francais]{babel}
\newtheorem{theorem}{Theorem}

\newtheorem{prop}{Proposition}

\newcommand{\R}{ {\mathbb R} }
\newcommand{\N}{{\mathbb N}}

\makeatletter
\@addtoreset{equation}{section}
\makeatother

\newcommand{\cqfd}{{\unskip\kern 6pt\penalty 500
\raise -2pt\hbox{\vrule\vbox to 6pt{\hrule width 6pt
\vfill\hrule}\vrule}\par}}

\begin{document}
\title{Convergence to equilibrium in competitive Lotka-Volterra equations}

\date{}
\author{Nicolas Champagnat$^{1}$, Pierre-Emmanuel Jabin$^{1,2}$,
  Ga\"el Raoul$^{3}$}

\footnotetext[1]{TOSCA project-team, INRIA Sophia Antipolis --
  M\'editerran\'ee, 2004 rte des Lucioles, BP.\ 93, 06902 Sophia
  Antipolis Cedex, France, \\
E-mail: \texttt{Nicolas.Champagnat@sophia.inria.fr}}

\footnotetext[2]{Laboratoire J.-A. Dieudonn\'e, Universit\'e de Nice --
  Sophia Antipolis, Parc Valrose, 06108 Nice Cedex 02, France, E-mail:
  \texttt{jabin@unice.fr}}
\footnotetext[3]{DAMTP, CMS, University of Cambridge, Wilberforce Road, 
Cambridge CB3 0WA, United Kingdom, Email: \texttt{g.raoul@damtp.cam.ac.uk}}
\maketitle{}
\selectlanguage{english}
\noindent{\bf Abstract.} We study a generalized system of ODE's
modeling a finite number of biological populations in a competitive
interaction. We adapt the techniques in \cite{JR} and
\cite{CJ} to prove the convergence to a unique stable equilibrium. 

\selectlanguage{french}
\noindent{\bf R\'esum\'e.} Nous \'etudions un syst\`eme
g\'en\'eralis\'e d'\'equations diff\'erentielles mod\'elisant un
nombre fini de populations biologiques en interaction
{com\-p\'e\-ti\-tive}. En adaptant les techniques de \cite{JR} et \cite{CJ},
nous prouvons la convergence vers un unique \'equilibre stable.

\bigskip

\noindent{\bf Version fran\c caise abr\'eg\'ee.}

Nous \'etudions le comportement en temps grand de mod\`eles de
dynamique de populations. On consid\`ere un nombre fini de
sous-populations, correspondant chacune \`a un trait ou type diff\'erent. Ces
populations interagissent entre elles de fa\c con comp\'etitive. En
notant $n_i(t)$ l'effectif de la sous-population num\'ero $i$, un des
mod\`eles les plus classiques est le syst\`eme de Lotka-Volterra comp\'etitif
\[
\frac{d}{dt} n_i=(r_i-\sum_j b_{ij}\,n_j)\,n_i, \quad i=1\ldots N,
\] 
o\`u $b_{ij}\geq 0$. On se place ici dans le cadre plus g\'en\'eral du
syst\`eme 
\[
\frac{d}{dt} n_i(t)=\left[
r_i-\int_{\Omega} K_i(\alpha)\;L\left(\sum_j B_j(\alpha)\,
n_j(t)\right)\,dP(\alpha)\right]\;n_i(t),\quad i=1\ldots N
\]
avec $(\Omega,P)$ un espace  mesurable. Ce syst\`eme peut
s'interpr\^eter comme un mod\`ele avec ressources g\'en\'eralis\'ees.

En utilisant les techniques d\'evelopp\'ees dans \cite{JR} pour une
version continue du premier mod\`ele, et dans \cite{CJ}, on peut
facilement montrer
\medskip

\noindent{\bf Th\'eor\`eme} {\em
Supposons que $L$ est une fonction $C^1$ sur $\R$, positive sur $\R_+$, que
$K$ et $B$ sont des fonctions positives appartenant \`a 
$L^\infty(dP(\alpha))\cap L^1(dP(\alpha))$ et que\\
\noindent $(i)$ (Comp\'etition stricte) 
$L$ est strictement croissante et pour tout $1\leq i\leq n$,
$r_i<\int_\Omega K_i(\alpha)L(\infty)dP(\alpha)$ 
o\`u $L(\infty):=\lim_{x\rightarrow+\infty}L(x)\in(0,+\infty]$.\\
\noindent $(ii)$ (Sym\'etrie) Il existe $C_i> 0$
tq $B_i(\alpha)=C_i\,K_i(\alpha)$ 
\\  
\noindent $(iii)$ (Non extinction)
 Pour tout $i$, $r_i>\int_{\Omega} K_i(d\alpha)L(0)\,dP(\alpha)$
\\
\noindent $(iv)$ (Non d\'eg\'en\'erescence) 
Pour $I\subset \{1\ldots N\}$, soit $\R^I$ l'ensemble des $n\in\R^N$
tels que $n_i=0$ pour 
tout $i\not\in I$. Pour tout $I\subset \{1\ldots N\}$ il y a au plus
un $n\in \R^{I}$ tq 
\[
r_i-\int_{\Omega} K_i(\alpha)\;L\left(\sum_{j=1}^N B_j(\alpha)\,
n_j\right)\,dP(\alpha)=0,\quad \forall i\in I.
\]
Alors $\exists!\;\bar n=(\bar n_1,\ldots, \bar n_N)\in
\R_+^N\setminus \{0\}$, tel que pour toute solution $n(t)=( n_1,\ldots,
n_N)$ du mod\`ele g\'en\'eralis\'e avec une donn\'ee initiale $n_i(0)>0$
$\forall i$, on a
\[
n(t)\longrightarrow \bar n, \ quand\ t\rightarrow +\infty.
\]
} 

\medskip
\noindent En particulier ce r\'esultat implique
\medskip

\noindent {\bf Proposition} {\em Supposons que $r_i>0$ pour tout $i$
  et que
la matrice $b_{ij}$ v\'erifie
\[
\exists C\in (\R_+^*)^N\ \mbox{tq}\ C_i\,
b_{ij}=b_{ji}\,C_j,\ \mbox{et}\ \sum_{ij}u_i\,u_j\,b_{ij}\,C_i>0\quad
\forall u\in\R^N\setminus \{0\},
\]
Alors $\exists!\;\bar n=(\bar n_1,\ldots, \bar n_N)\in
\R_+^N\setminus\{0\}$ tel que pour toute solution $n(t)=( n_1,\ldots,
n_N)$ du premier mod\`ele avec donn\'ee initiale $n_i(0)>0$ $\forall i$, 
\[
n(t)\longrightarrow \bar n, \ quand\ t\rightarrow +\infty.
\]
} 
\selectlanguage{english}
\section{Introduction}
We study the long time behaviour of models of population dynamics. We
consider a finite number of 
subpopulations whose dynamics is governed by a system of competitive
ODEs (in the sense of Hirsch, see 
e.g.~\cite{hirsch-88}). We denote by $n_i(t)$, $i=1\ldots N$, the
number of individuals 
of the subpopulation $i$.

The most classical models are competitive Lotka-Volterra equations
\begin{equation}
\frac{d}{dt} n_i=(r_i-\sum_j b_{ij}\,n_j)\,n_i, \quad i=1\ldots N,\label{Eql}
\end{equation}
where $b_{ij}\geq 0$, and the models with a finite number of resources
\begin{equation}
\frac{d}{dt} n_i=(-d_i+\sum_{k=1}^K I_k\,\eta_{ki})\,n_i, \label{Eqr}
\end{equation}
where $\eta_{ki}\geq 0$ and the $I_k$ are given by the
Holling II functional response 
\begin{equation*}
I_k=\frac{I_k^0}{1+\sum_{i=1}^N \eta_{ki}\,n_i}.
\end{equation*}
This type of system appears in biology when one studies the dynamics
of a system of interacting species (see 
\cite{HS,Gop,Kri}). It also appears in Trait Substitution Sequence
models, where one considers a population 
structured by a continuous phenotype (see equation \eqref{Eqlcont},
\eqref{continuous} on this matter), where 
only a small number of traits are present (see
\cite{MGMJH,Cha}). These models have been used to develop the 
theory of Adaptative Dynamics (see \cite{MGMJH,Cha,Di}). 

Previous asymptotic studies on this type of equations concern either
very general properties (the existence of 
a carrying simplex~\cite{hirsch-88}), or precise results but only for
low dimensional systems ($N\leq 
3$~\cite{zeeman-93}), under strong assumptions of the coefficients
(for instance, the matrix $(b_{ij})$ is 
supposed to be diagonal dominant, see \cite{HS}), or only on local
properties (the equilibrium population is 
locally stable, or populations $n_i$ do not vanish).

Note that both equations \eqref{Eql} and \eqref{Eqr} may be
interpreted as discrete versions of 
continuous models. To each subpopulation corresponds a phenotypic
trait $x_i\in \R^d$, and then posing
\begin{equation}
n(t,x)=\sum_{i=1}^N n_i(t) \delta_{x_i},\label{dirac}
\end{equation}
one finds that Eq. \eqref{Eql} for instance is equivalent to
\begin{equation}
\partial_t n(t,x)=(r(x)-\int_{\R^d} b(x,y)\,n(t,dy))\,n(t,x),\label{Eqlcont}
\end{equation}
with $r_i=r(x_i)$ and $b_{ij}=b(x_i,x_j)$.

The long time behaviour of the continuous model \eqref{Eqlcont}
(with bounded initial 
data instead of Dirac masses) was studied in \cite{JR}. For a
symmetric $b$ defining a positive operator,
the convergence to the unique stable equilibrium was proved. For the
case with resources, the result is essentially contained in \cite{CJ},
which generalizes the derivation of \cite{DJMP}.

The study of the discrete or continuous models corresponds to slighty
different biological questions; in the 
continuous case, it is for instance connected to the issue of
speciation, or how from a continuum of traits a 
few well separated ones (the ``species'') are selected; in the
discrete case, one is rather concerned about 
survival or extinction of each subpopulations. From a rigorous
mathematical point of view, a result in the 
continuous case does not imply anything for the discrete one.  However
it is easy to apply the techniques developed in \cite{JR}
and \cite{CJ} to the discrete models; 
that is our aim.

First of all, 
we consider the very general equation
\begin{equation}
\frac{d}{dt} n_i(t)=\left[
r_i-\int_{\Omega} K_i(\alpha)\;L\left(\sum_j B_j(\alpha)\,
n_j(t)\right)\,dP(\alpha)\right]\;n_i(t),\quad i=1\ldots N,\label{discrete}
\end{equation}
with $(\Omega,P)$ any measurable space,  or in the continuous case
\begin{equation}
\partial_t n(t,x)=\left[
r(x)-\int_{\Omega} K(x,\alpha)\;L\left(\int_{\R^d} B(y,\alpha)\,
n(t,dy)\right)\,dP(\alpha)\right]\;n(t,x).\label{continuous}
\end{equation}
We prove the following
\begin{theorem}
  Assume that $L$ is $C^1$ on $\R$ and non negative on $\R_+$, that
  $K$ and $B$ are non negative, in $L^\infty(dP(\alpha))\cap
  L^1(dP(\alpha))$ and that\\ 
  \noindent $(i)$ (Strict competition) $L$ is strictly increasing and
  $r_i<\int_\Omega 
  K_i(\alpha)L(\infty)\,dP(\alpha)$ for all $1\leq i\leq N$, where
  $L(\infty):=\lim_{x\rightarrow+\infty}L(x)\in(0,+\infty]$.\\
  \noindent $(ii)$ (Symmetry) There exists $C_i>0$
  s.t. $B_i(\alpha)=C_i\,K_i(\alpha)$ 
  \\  
  \noindent $(iii)$ (Non extinction)
  For any $i$, $r_i>\int_{\Omega} K_i(\alpha)L(0)\,dP(\alpha)$
  \\
  \noindent $(iv)$ (Non degeneracy) 
  For any subset $I\subset \{1\ldots N\}$, let $\R^I$ be the set of
  $n\in\R^N$ s.t.\ $n_i=0$ for all $i\not\in 
  I$. For all $I\subset\{1\ldots N\}$, there exists
  at most one $n\in \R^{I}$ s.t.\
  \begin{equation}
    \label{eq:cond4}
  r_i-\int_{\Omega} K_i(\alpha)\;L\left(\sum_{j=1}^N B_j(\alpha)\,
    n_j\right)\,dP(\alpha)=0,\quad \forall i\in I.
  \end{equation}
  Then there exists a unique $\bar n=(\bar n_1,\ldots, \bar n_N)\in
  \R_+^N$ with $\bar n\neq 0$, s.t. for any solution $n(t)=( n_1,\ldots,
  n_N)$ to \eqref{discrete} with initial data $n_i(0)>0$ for any $i$, 
  \[
  n(t)\longrightarrow \bar n, \ as\ t\rightarrow +\infty.
  \]
  \label{mainthm}
\end{theorem}
Note that, if one had $r_i>\int K_i(\alpha)L(\infty)\,dP(\alpha)$ for
some $i$, then $n_i(t)\rightarrow+\infty$ if
$n_i(0)>0$. Assumption~$(i)$ hence ensures the non-explosion of the system.

Eq.~\eqref{discrete} could be directly derived from simple biological
considerations. It assumes that the reproduction rate of a population
of type $i$ (or with trait $x_i$) is the difference between a fixed rate
depending only on the trait and a competitive interaction with the
other populations, resulting from the interaction with
the environment.  The state of each component of this environment
(indicated by different values of $\alpha$) is given by the sum
\[
\sum_j B_j(\alpha)\,n_j(t). 
\] 
Each such component has some independent effect on the reproduction. 
To get the total reproduction rate
one sums over those.

Eq.~\eqref{discrete} is hence an obvious generalization, with a
possibly infinite number of resources, of the 
model \eqref{Eqr}.  It also contains the Lotka-Volterra
system~\eqref{Eql}. In this case, 
Theorem~\ref{mainthm} gives
\begin{prop} Assume that $r_i>0$ for all $i$ and that
the matrix $b_{ij}$ satisfies
\begin{equation}
\exists C\in (\R_+^*)^N\ \mbox{s.t.}\ C_i\,
b_{ij}=b_{ji}\,C_j,\ \mbox{and}\ \sum_{ij}u_i\,u_j\,b_{ij}\,C_i>0\quad
\forall u\in\R^N\setminus \{0\}, \label{hypb}
\end{equation}
then there exists a unique $\bar n=(\bar n_1,\ldots, \bar n_N)\in
\R_+^N$ with $\bar n\neq 0$, s.t.\ for any solution $n(t)=( n_1,\ldots,
n_N)$ to \eqref{Eql} with initial data $n_i(0)>0$ for any $i$, 
\[
n(t)\longrightarrow \bar n, \ as\ t\rightarrow +\infty.
\]
\label{propb}
\end{prop}


  This result shows that, in Lotka-Volterra systems which are symmetric in the sense of~(\ref{hypb}), the
  competition between a mutant trait and a resident population leads to a unique stationary state,
  regardless of the initial population state. This is precisely the assumption needed in~\cite{Cha} to apply a
  limit of large population and rare mutations to an individual-based
  model.
  In particular, Thm.~2.7 of~\cite{Cha} applies to symmetric competitive Lotka-Volterra systems.

\bigskip
\noindent {\bf Proof of Prop.~\ref{propb}.}
Define the matrix $m_{ij}=C_i\,b_{ij}$. Note that $m$ is symmetric and
positive definite. Hence there exists an orthonormal basis of
eigenvectors $U^i$,
$i=1\dots N$, and corresponding eigenvalues $\lambda_i>0$. 

Then put $L=\text{Id}$, 
$\Omega=\{1,...,N\}$, $P=\frac{1}{N}\sum_{i=1}^N \delta_i$, 
$B_j(\alpha)=\sqrt{\lambda_\alpha}
U^\alpha_j$, $K_i(\alpha)=C_i^{-1}\,\sqrt{\lambda_\alpha}U^\alpha_i$
and note that
\[\begin{split}
\sum_{j=1}^N b_{ij}\,n_j&=\frac{1}{C_i}\sum_{j=1}^N
m_{ij}\,n_j=\frac{1}{C_i}\left[Mn\right]_i\\
&=\frac{1}{C_i}\left[M\left(\sum_\alpha U^\alpha\langle
  U^\alpha,n\rangle\right)\right]_i= 
\frac{1}{C_i}\sum_{\alpha=1}^N \lambda_\alpha\,U^\alpha_i \left(\sum_{j=1}^N
U^\alpha_j\,n_j\right)\\
&=\int_\Omega K_i(\alpha) L\left(\sum_{j=1}^N
B_j(\alpha)\,n_j\right)\,dP(\alpha). 
\end{split}
\]
Therefore Eq.~\eqref{discrete} indeed yields \eqref{Eql} in that
particular case.

Conditions $(i)$ and $(ii)$ of Theorem~\ref{mainthm} are obviously
satisfied. Conditions $(iii)$ holds since $r_i>0$ for all $i$ and
$L(0)=0$. As for condition $(iv)$, assume that for a subset $I$ one
has two vectors $n_j^\gamma$, $\gamma=1,2$, s.t. $n_i^\gamma=0$ for
$i\not\in I$ and
\[
r_i=\int_\Omega K_i(\alpha) L\left(\sum_{j=1}^N
B_j(\alpha)\,n_j^\gamma\right)\,dP(\alpha)=\sum_{j=1}^N b_{ij}
n_j^\gamma\quad \forall 
i\in I.
\]
Put $\delta n=n^1-n^2$ and simply note that
$\delta n_i \sum_{j=1}^N b_{ij}\delta n_j =0$ pour $ i=1\dots N$.

This means that $\delta n=0$ and proves $(iv)$ since
$\sum_{i,j=1}^N C_i b_{ij}\,\delta n_i\,\delta n_j=0$.
Hence the proposition is implied by Theorem~\ref{mainthm}.
\cqfd
\medskip

Note that the same argument works in the continuous case and
Eq.~\eqref{Eqlcont} is a particular case of~\eqref{continuous} for
$x\in O$ a bounded domain. The condition on $b$ is 
\[
C(x)\,b(x,y)=C(y)\,b(y,x),\quad \int_{O^2}
C(x)\,b(x,y)\,n(x)\,n(y)\,dx\,dy>0\quad \forall n\neq 0.
\]
One still puts $L(\xi)=\xi$. 
Notice that $C(x)\,b(x,y)$ defines a
 compact, selfadjoint and positive operator on $L^2(O)$. 
Diagonalizing the operator,
one gets 
\[
C(x)\,b(x,y)=\sum_{\alpha} \lambda_i f_\alpha(x)\,f_\alpha(y),
\]
with $\lambda_\alpha>0$ the eigenvalues, tending to $+\infty$ and
$f_\alpha$ the corresponding normalized eigenvector. It is hence enough to
take $\Omega=\N$ and
$
K(x,\alpha)=\sqrt{\lambda_\alpha}\,f_\alpha(x). 
$

In the particular case where $b(x,y)=b(x-y)$ on the whole $\R^d$
and $C=1$, by Fourier transform,
the condition on $b$  means that $\hat b>0$.
One then takes $\Omega=\R^d$ and
\[
K(x,\alpha)=(\cos(\alpha\cdot x)\,\sqrt{\hat
  b(\alpha)},\ \sin(\alpha\cdot x)\,\sqrt{\hat b(\alpha)} ).
\]
%
\section{Proof of Theorem \ref{mainthm}}
The proof is based on the study of the following Lyapunov functional
\begin{equation}
F(n)=\int_\Omega H\left(\sum_{j=1}^N B_j(\alpha)
n_j\right)\;dP(\alpha)-\sum_{i=1}^N C_i\, r_i\, n_i,\label{F}
\end{equation}
where $H$ is an antiderivative of $L$ and hence strictly convex.
\subsection{$F$ is a strict Lyapunov functional}\label{sub1}
Let $n$ be a solution to \eqref{discrete}. Then by a direct
computation 
\begin{align}
\frac{d}{dt} F(n(t))
&=-\sum_{i=1}^N C_i\, n_i\,\left[
\int_\Omega K_i(\alpha)\,L\left(\sum_{j=1}^N B_j(\alpha)
n_j\right)\;dP(\alpha)-\,r_i
\right]^2. \label{eq:F-Lyap}
\end{align}
Therefore $F(n(t))$ is non increasing and its derivative in time vanishes only on stationary solutions
to~\eqref{discrete}, i.e.\ $F$ is a strict Lyapunov functional for the system~(\ref{discrete}).  

Thanks to condition~$(i)$,
\begin{align}
  \frac{\partial F}{\partial n_i} & 
  \geq C_i\left(\int_\Omega K_i(\alpha)L\Big(\frac{K_i(\alpha)}{C_i}
  n_i\Big)\,dP(\alpha)-r_i\right) 
  \geq a>0 \notag
\end{align}
if $n_i$ is large enough. Therefore, there is a constant $a'>0$ s.t.\ $\nabla F(n)\cdot n\geq a'\|n\|$ if $\|n\|$ is large enough. This implies that $F(n)\rightarrow+\infty$ when $\|n\|\rightarrow+\infty$, and entails that $n(t)$ is uniformly bounded.

Let $n\in\R_+^N$ be a steady-state of \eqref{discrete} and let $I$ be the set of $i$ s.t.\ $n_i>0$. Then, for
any $i\in I$ one needs to have
\[
\int_\Omega
K_i(\alpha)\, 
L\left(\sum_{j=1}^N B_j(\alpha) 
n_j\right)\;dP(\alpha)=r_i.
\] 
By condition $(iv)$ there is at most one such solution for every $I$, and there are only a finite number of
possible $I$, $F$ has then a finite number of steady-states.

Classical Lyapunov functionnals' techniques 
then entail that the solution $n(t)$ to~(\ref{discrete}) converges to a steady-state $\tilde
n$ for any initial condition $n(0)$.

\subsection{The functional $F$ is convex}

Compute
\begin{equation}
\label{eq:der2-F}
\frac{\partial^2 F}{\partial n_i\partial n_k}=\int_\Omega
B_i(\alpha)\,B_k(\alpha) 
L'\left(\sum_{j=1}^N B_j(\alpha) 
n_j\right)\;dP(\alpha).
\end{equation}
Hence as $L$ is increasing
\begin{equation}
\label{eq:unif-ellipt}
\sum_{i,k} \frac{\partial^2 F}{\partial n_i\partial n_k} \xi_i\,\xi_k=
\int_\Omega
(\sum_i \xi_i B_i(\alpha))^2\, 
L'\left(\sum_{j=1}^N B_j(\alpha) 
n_j\right)\;dP(\alpha)\geq 0.
\end{equation}
$F$ is therefore convex and any local minimum on $\R_+^N$ is global.

Since~(\ref{discrete}) has a finite number of stationary solutions, this clearly implies that $F$ admits a
unique global minimizer $\bar{n}$. Otherwise, $F$ would reach its minimum on the whole segment linking two
distinct minimizers.

The object of the next subsection is to prove that $\bar{n}$ satisfies a stronger property: this is the unique
ESS of the system.

\subsection{Uniqueness of the ESS}\label{sub2}

%
%


Any local minimizer  $n\in\mathbb R_+^N$ of the functional $F$  necessarily satisfies
\begin{equation}
\begin{split}
&\int_\Omega
K_i(\alpha)\, 
L\left(\sum_{j=1}^N B_j(\alpha) 
n_j\right)\;dP(\alpha)=r_i,\quad \forall i\ s.t.\ n_i>0,\\
&\int_\Omega
K_i(\alpha)\, 
L\left(\sum_{j=1}^N B_j(\alpha) 
n_j\right)\;dP(\alpha)\geq r_i,\quad \forall i\ s.t.\ n_i=0.
\end{split}\label{ESS}
\end{equation}
This condition corresponds to the usual definition of an Evolutionarily Stable Strategy in adaptive dynamics
(see for instance \cite{Di}).  It turns out that there exists at most one ESS, $\bar{n}$. Hence being an ESS
is a necessary and sufficient condition to be the global minimizer of $F$.

Indeed take two $n^\gamma\in \R^N_+$, $\gamma=1,2$ satisfying \eqref{ESS} and
compute
\[\begin{split}
 0&\geq \sum_i C_i\,n^1_i \left(r_i-\int_\Omega
K_i(\alpha)\, 
L\left(\sum_{j=1}^N B_j(\alpha) 
n_j^2\right)\;dP(\alpha)\right)\\
&\qquad+\sum_i C_i\,n^2_i \left(r_i-\int_\Omega
K_i(\alpha)\, 
L\left(\sum_{j=1}^N B_j(\alpha) 
n_j^1\right)\;dP(\alpha)\right). \\
\end{split}\]
This last quantity is equal to (thanks to \eqref{ESS}) 
\[
\begin{split}
& \sum_i C_i\,(n^1_i-n^2_i) \left(r_i-\int_\Omega
K_i(\alpha)\, 
L\left(\sum_{j=1}^N B_j(\alpha) 
n_j^2\right)\;dP(\alpha)\right)\\
&\qquad+\sum_i C_i\,(n^2_i-n^1_i) \left(r_i-\int_\Omega
K_i(\alpha)\, 
L\left(\sum_{j=1}^N B_j(\alpha) 
n_j^1\right)\;dP(\alpha)\right) \\
\end{split}\]
and to
\[\begin{split}
&+ \int_\Omega \left(\sum_j B_j(\alpha)\,n_j^1-\sum_j B_j(\alpha)
n_j^2\right)\\
&\qquad\left(L\left(\sum_{j=1}^N B_j(\alpha) 
n_j^1\right)-L\left(\sum_{j=1}^N B_j(\alpha) 
n_j^1\right)\right)\;dP(\alpha).
\end{split}\] 
As $L$ is strictly increasing, this implies that for $P\ a.e.$
$\alpha$,
$\sum_{i=1}^N B_i(\alpha)\,(n^1_i-n^2_i)=0$ and by   
$(iv)$, it means that $n^1=n^2$.
\subsection{Conclusion of the proof of Thm.~\ref{mainthm}}\label{sub3}

Assume that $n_i(0)>0$ for all $1\leq i\leq N$. We know from Subsection \ref{sub1} that $n(t)$ converges to a
steady-state $\tilde n$ when $t\to\infty$.

If $\tilde n$ does not satisfy \eqref{ESS}, there exists $i\in \{1,\dots,N\}$ such that
\[
\lambda_i:=r_i-\int_\Omega
K_i(\alpha)\, 
L\left(\sum_{j=1}^N B_j(\alpha) 
\tilde n_j\right)\;dP(\alpha)> 0.
\]
Since $n_i(0)>0$, $n_i>0$ at all times, and the linearized equation around $\tilde n$ shows that $n$ cannot converge to $\tilde n$:
\begin{eqnarray*}
  \label{eq:}
  \frac d{dt}(n-\tilde n)_i&=&\left(\lambda_i+O(\|n-\tilde n\|)\right)\,(n-\tilde n)_i\\
&\geq&\frac{\lambda_i}2\,(n-\tilde n)_i,
\end{eqnarray*}
provided that $\|n-\tilde n\|$ is small enough.

Therefore, $\tilde{n}=\bar{n}$, and the proof of Thm.~\ref{mainthm} is completed.

\bigskip

\noindent{\bf Acknowledgments:} The first author is grateful to Michel Bena\"im for useful discussions on the
dynamical systems context of the problem. GR has been supported by Award No. KUK-I1-007-43 of Peter A. Markowich, 
made by King Abdullah University of Science and Technology (KAUST).

\end{document}